\documentclass{article}
\title{Configuration Spaces and ${\mathbb R}^n$}
\author{Jeffrey L. Caruso}
\date{September 16, 2004}

\usepackage{amssymb}
\usepackage{exscale}
\usepackage{stmaryrd}
\usepackage{mathrsfs}
\usepackage[all]{xy}

\newcommand{\rar}{\rightarrow}

\newcommand{\CY}[1]{\mathscr{C}(Y)_{#1}}
\newcommand{\SJ}{\mathcal{S}_j}
\newcommand{\RxY}{\mathbb{R}\times Y}
\newcommand{\RnxY}{\mathbb{R}^n\times Y}
\newcommand{\CRY}[1]{\mathscr{C}(\mathbb{R}\times Y)_{#1}}
\newcommand{\ConeY}[1]{\mathscr{C}_1(Y)_{#1}}

\newcommand{\CbarYX}{\bar{C}_1(Y,X)}
\newcommand{\EYX}{E_1(Y,X)}
\newcommand{\CYSX}{C(Y,\Sigma X)}
\newcommand{\OnCYSnX}{\Omega^n C(Y,\Sigma^n X)}
\newcommand{\AR}{\mathbb{R}}
\newcommand{\half}{\frac{1}{2}}
\newcommand{\threeq}{\frac{3}{4}}
\newcommand{\qtr}{\frac{1}{4}}

\newtheorem{thm}{Theorem}
\newtheorem{cor}{Corollary}
\newtheorem{defn}{Definition}

\begin{document}

\maketitle

This note attempts to make clear the relation between configurations of 
points in a space $Y$ and those in its Cartesian product with the reals. It 
turns out to be a very simple relation whose proof uses nothing new.

Let $Y$ be an unbased space. Denote by $Y^j$ the $j$-fold Cartesian product 
of $Y$ with itself. For present purposes we consider the circle $S^1$ to 
be the quotient of the unit interval $[0,1]/\{0,1\}$. If $X$ is a based space 
then $\Sigma X$ is defined to be $X\wedge S^1$ and $\Omega X$ is defined 
to be the loop space of $X$, that is, the space of based maps from $S^1$ 
to $X$.

\begin{defn}
Define $\CY{j}$ to be the subspace of $Y^j$ consisting of $j$-tuples of 
distinct points in $Y$.  If $\nu$ is an injective function from 
$\{1,\ldots,i\}$ to $\{1,\ldots,j\}$ then define $\nu^*: \CY{j} \rar \CY{i}$ 
by sending $(y_1,\ldots,y_i)$ to $(y_{\nu(1)},\ldots,y_{\nu(i)})$. If $X$ is 
a nondegenerately based space, define $\nu_* : X^i \rar X^j$ sending 
$(x_1,\ldots,x_i)$ to $(x_1',\ldots,x_j')$ where $x_k' = x_l$ if $k = \nu(l)$
and $x_k'$ is the basepoint if $k$ is not in the image of\ $\nu$.
\end{defn}

Note that these maps are compatible with composition; i.e. 
$(\nu\circ\mu)_* = \nu_*\circ\mu_*$ and $(\nu\circ\mu)^* = \mu^*\circ\nu^*$. 
In particular, the maps $\nu^*$ define a free action of the $j$-fold 
symmetric group ${\cal S}_j$ on $\CY{j}$.

The spaces $\CY{j}$ and the maps $\nu^*$ define a {\it coefficient system} 
in the sense of \cite{cmt}, 
and we define an equivalence relation $\sim$ on 
$\coprod_j \CY{j} \times X^j$ generated by 
$(\nu^*(\vec y),\vec x) \sim (\vec y,\nu_*(\vec x))$. Define
$$C(Y,X) = \left(\coprod_j \CY{j} \times X^j\right) \Big/ {\sim}.$$

In their recent paper \cite{ct}, Cohen and Taylor deal 
with the space $C(\RxY,X)$. Recall that a 
{\it weak metric space} is a space $Y$ together with a continuous 
function $d:Y\times Y\rar [0,\infty)$ such that $d^{-1}(0)$ is the 
diagonal in $Y\times Y$. The main result of this note is:

\begin{thm}
Let $Y$ be a weak metric space and $X$ a nondegenerately based space.
There is a space $C_1(Y,X)$ and a pair of maps 
$$ C(\RxY,X) {\buildrel {\phi} \over \longleftarrow} C_1(Y,X)
{\buildrel {\alpha} \over \longrightarrow} \Omega\CYSX $$
such that:
\begin{enumerate}
\item $C_1(-,-)$ is functorial with respect to based maps in the second 
variable and injective maps in the first variable, and $\phi$ and $\alpha$ 
are natural;
\item $\phi$ is a homotopy equivalence; and
\item $\alpha$ is a weak homotopy equivalence if $X$ is path-connected.
\end{enumerate}
\end{thm}

The proof uses the methods from \cite{may} and \cite{cmt}.
The space $C_1(Y,X)$ is another 
space derived from a ``coefficient system.'' Let $\ConeY{j}$ be 
the subspace of $(\AR\times\AR\times Y)^j$ consisting of $j$-tuples
of triples $\big((a_1,b_1,y_1),\ldots,(a_j,b_j,y_j)\big)$ such that 
for all $i$, $a_i < b_i$ and for all $k \ne l$, $y_k = y_l$ implies 
$b_k \le a_l$ or $b_l \le a_k$.

We can define a coefficient system structure $\nu^*$, $\nu_*$ on 
$\big\{\ConeY{j}\big\}_{j\ge 0}$ by acting on triples, and define
and $\sim$ on $\coprod_j \ConeY{j} \times X^j$ generated by 
$(\nu^*(\kappa),\vec x) \sim (\kappa,\nu_*(\vec x))$. The quotient
space $C_1(Y,X)$ can be thought of as consisting of configurations of 
line segments in $\RxY$ with disjoint interiors, labeled by points of $X$; 
a segment labeled by the basepoint drops out under the identification\ 
$\sim$. For compactness of notation, we will use 
$$\big(a_i,b_i,y_i\big)_{1\le i\le j}$$
as shorthand for $\big((a_1,b_1,y_1),\ldots,(a_j,b_j,y_j)\big) \in 
\ConeY{j}$, and  
$$\big[a_i,b_i,y_i,x_i\big]_{1\le i\le j}$$
for the image of $\big((a_1,b_1,y_1),\ldots,(a_j,b_j,y_j),
(x_1,\ldots,x_j)\big)$ in $C_1(Y,X)$. Similarly we will use the shorthand
$\big[y_i,x_i\big]_{1\le i\le j}$ for points of $C(Y,X)$.

There is an obvious map $\phi_j$ from $\ConeY{j}$ to $\CRY{j}$ 
taking each segment to its center-point. This map respects permutations 
and so induces a map $\phi$ from $C_1(Y,X)$ to $C(\RxY,X)$.

There is also a map $\bar{\phi}_j$ from $\CRY{j}$ to $\ConeY{j}$ 
which we define as follows. Use the weak metric $d$ on $Y$ to 
define $g:(\RxY)\times(\RxY)\rar[0,\infty)$ by setting 
$$g\big((a,y),(a',y')\big) = \half \left(
\frac{\vert a-a'\vert^2 + d(y,y')}{\vert a-a'\vert + d(y,y')+1}
\right)$$
so $g((a,y),(a',y'))\le\half\vert a-a'\vert$ if $y=y'$. Let $\kappa = 
((a_1,y_1),\ldots,(a_j,y_j)) \in \CRY{j}$ and define 
$$v(\kappa) = \min_{k\ne l}\big\{g((a_k,y_k),(a_l,y_l))\big\}.$$
It's clear that $v(\kappa)>0$ and that the intervals $[a_k - v(\kappa), 
a_k + v(\kappa)]$ and $[a_l - v(\kappa), a_l + v(\kappa)]$ do not overlap 
when $y_k=y_l$, so we can define 
$$\bar{\phi}_j(\kappa) = 
\big(a_i-v(\kappa),a_i+v(\kappa),y_i\big)_{1\le i\le j}.$$
These induce a map $\bar{\phi}: C(\RxY,X)\rar C_1(Y,X)$. Further, $\phi_j$ 
and $\bar{\phi}_j$ are easily seen to be inverse $\SJ$-equivariant homotopy
equivalences: $\phi_j\bar{\phi}_j$ is the identity of $\CRY{j}$, and there 
is a deformation from the identity of $\ConeY{j}$ to $\bar{\phi}_j\phi_j$
by linearly scaling the intervals around their centers. 
So by Lemma 2.7(ii) of \cite{cmt}, $\phi$ is a homotopy equivalence.

Next we need to define $\alpha$. For the purposes of this section it
is more convenient to work with a homeomorphic copy of $C_1(Y,X)$. Let
$$\CbarYX = \left\{\big[a_i,b_i,y_i,x_i\big]_{1\le i\le j}\in C_1(Y,X)
\big| 0<a_i<b_i<1\hbox{ for all $i$}\right\}$$
This subspace is clearly homeomorphic to $C_1(Y,X)$ via the homeomorphism
of the reals $\AR$ with the open interval $(0,1)$.
Let $w=[a_i,b_i,y_i,x_i]_{1\le i\le j}$ be a point of $\CbarYX$. 
For a given $t$, define
$$\alpha(w)(t) = 
\big[y_i,[x_i,s_i]\big]_{1\le i\le j\hbox{ and }a_i\le t\le b_i},$$
where $s_i = (t-a_i)/(b_i-a_i)$. For a given $t$ and for each $i$ 
satisfying $a_i\le t\le b_i$, we observe that $0\le s_i\le 1$
and the points 
$\left\{y_i\big| 1\le i\le j\hbox{ and }a_i\le t\le b_i\right\}$ are 
distinct; also $\alpha(w)(0) = \alpha(w)(1)$ is the basepoint $*$ of 
$C(Y,\Sigma X)$. Thus $\alpha(w)$ is a well-defined loop in $\CYSX$.

To show $\alpha$ is a weak equivalence, we use the same idea as \cite{may}, 
namely to fit it into a comparison of quasifibration sequences. Define 
$\EYX$ to be the quotient space of $\CbarYX\times [0,1]$ where we identify 
$\left(\big[a_i,b_i,y_i,x_i\big]_{1\le i\le j},s\right)$ and 
$\left(\big[a'_i,b'_i,y'_i,x'_i\big]_{1\le i\le j+k},s\right)$
if $(a_i,b_i,y_i,x_i)=(a'_i,b'_i,y'_i,x'_i)$ for $1\le i\le j$ and
$a_{j+1},\ldots,a_{j+k} \ge s$. Note that all points of the form 
$(w,0)$ are identified with the basepoint $(*,0)$ of $\EYX$, so $\EYX$ is
contractible.

Define a map $\bar{\alpha}$ from $\EYX$ to the path space $P\CYSX$ by 
$$\bar{\alpha}(w,s)(t)=\cases{\alpha(w)(t),&\hbox{if $t\le s$, and}\cr
\noalign{\vskip2pt}\alpha(w)(s),&\hbox{if $t\ge s$.}\cr}$$

Defining $\iota:\CbarYX\rar \EYX$ by $\iota(w) = (w,1)$ and $q:\EYX\rar 
\CYSX$ by $q(w,s) = \bar{\alpha}(w,s)(1)$, we have the following commutative
diagram
$$\xymatrix{
\CbarYX \ar[r]^\alpha \ar[d]_{\iota} & \Omega\CYSX \ar@{^{(}->}[d] \\
\EYX \ar[r]^{\bar{\alpha}} \ar[d]_{q} & P\CYSX \ar[d]^{p_1} \\
\CYSX \ar@{=}[r] & \CYSX }$$
where $p_1$ is projection on the endpoint. Thus by comparison of the long 
exact sequences of homotopy groups, it is enough to show that $q$ is a 
{\it quasifibration}, that is, a map $q:E\rar B$ such that for all $b\in B$ 
the canonical map from $q^{-1}(b)$ to the homotopy fiber of $q$ over $b$ is
a weak homotopy equivalence.

Recall from \cite{dt} the Dold-Thom criterion for a map over a filtered base 
space to be a quasifibration. Let $B$ be a space with closed 
subspaces
$$F_{0}B\subseteq F_{1}B\subseteq \ldots F_{j}B\subseteq \ldots \subseteq B$$
and $B = \bigcup_{j\ge 0}F_{j}$, and let $q:E\rar B$ be a map. 
A subspace $V\subseteq B$ is called {\it distinguished} if the 
restriction $q:q^{-1}(V)\rar V$ is a quasifibration. Then

\begin{thm} (Dold and Thom) $B$ is distinguished provided that
\begin{enumerate}
\item $F_{0}B$ is distinguished, and for each $j>0$ every open subset of 
$F_{j}B\smallsetminus F_{j-1}B$ is distinguished, and
\item for each $j>0$ there is a homotopy $h_t:U\rar U$ of a neighborhood $U$ 
of $F_{j-1}B$ in $F_{j}B$, and a homotopy $H_t: q^{-1}(U)\rar q^{-1}(U)$ 
such that:
\begin{enumerate}
\item $h_0$ is the identity map of $U$, $h_1(U)\subseteq F_{j-1}B$,
and for all $t$,\\$h_t(F_{j-1}B)\subseteq F_{j-1}B$,
\item $H_0$ is the identity map of $q^{-1}(U)$
and for all $t$, $qH_t=h_t q$, and
\item for all $z\in U$, the map 
$H_1:q^{-1}(z)\rar q^{-1}(h_1(z))$
is a homotopy equivalence.
\end{enumerate}
\end{enumerate}
\end{thm}

Here we give $\CYSX$ the filtration of \cite{may}, that is $F_j\CYSX$ is 
defined to be the image of $(\coprod_{0\le k\le j} \CY{k} \times 
(\Sigma X)^k)$. This has the property that $F_0\CYSX$ consists of just
the basepoint *, and $F_j\CYSX\smallsetminus F_{j-1}\CYSX$ is 
homeomorphic to the image of
$\ConeY{j}\times \big(X\smallsetminus \{*\}\times (0,1)\big)^j$.

We define some maps on $\CbarYX$ to help elucidate the proof. If 
$w=\big[a_i,b_i,y_i,x_i\big]_{1\le i\le j}$ and 
$w'=\big[a_i,b_i,y_i,x_i\big]_{j+1\le i\le j+k}$ are configurations 
in which for all $k\ne l$ the sets 
$$\left\{(t,y_i)\in \RxY \big| a_i < t < b_i\right\}$$
are pairwise disjoint, then let
 $w\cup w'=\big[a_i,b_i,y_i,x_i\big]_{1\le i\le j+k}$.
This is continuous on the subspace of $\CbarYX\times\CbarYX$ on which it is
defined.

If $s$ and $t$ are real numbers with $s<t$ and 
$w=\big[a_i,b_i,y_i,x_i\big]_{1\le i\le j}$, then define
$$shrink_{s,t}(w)=\big[s+(t-s)a_i,s+(t-s)b_i,y_i,x_i\big]_{1\le i\le j},$$
which linearly compresses a configuration of segments in $(0,1)\times Y$
into the slice $(s,t)\times Y$. Note that the composition 
$\mu:\CbarYX\times\CbarYX\rar\CbarYX$ defined by 
$$\mu(w,w')=shrink_{0,\half}(w)\cup shrink_{\half,1}(w')$$
defines an $H$-space structure on $\CbarYX$.

For an element $z=\big[y_i,[x_i,s_i]\big]_{1\le i\le j}\in 
F_j\CYSX\smallsetminus F_{j-1}\CYSX$ we define 
$$\lambda(z)=
\left[\half-\frac{s_i}{2},1-\frac{s_i}{2},y_i,x_i\right]_{1\le i\le j}.$$
This maps via $\alpha$ to a loop whose value is $z$ at $t=\half$,
and is well-defined and continuous on $F_j\CYSX\smallsetminus F_{j-1}\CYSX$.

For an element $w\in\CbarYX$ and $s\in [0,1]$, we can define a function
$$below_s(w)=\big[a_i,b_i,y_i,x_i\big]_{1\le i\le j\hbox{ and }b_i\le s},$$
the segments of $w$ contained in $[0,s]\times Y$. This is continuous
on $q^{-1}(F_j\CYSX\smallsetminus F_{j-1}\CYSX)$.

For a relatively open set $V\subseteq F_j\CYSX\smallsetminus F_{j-1}\CYSX$ 
define $\psi:\CbarYX\times V\rar q^{-1}(V)$ by 
$$\psi(w,z)=\left(shrink_{0,\half}(w)\cup shrink_{\half,1}(\lambda(z)),
\threeq\right).$$

If $(w,s)\in\EYX$ define $\bar{\psi}(w,s) = below_s(w).$ It follows that
there is a commutative diagram
$$\xymatrix{
\CbarYX\times V \ar@<-3pt>[rr]_{\psi} 
\ar[dr]_{p_2} &  & q^{-1}(V) \ar@<-3pt>[ll]_{(\bar{\psi},q)} \ar[dl]^q \\
& V } $$
The left map is projection on the second factor and so is the
simplest kind of quasifibration; thus the proof of part (1.) 
will be complete when we have shown that $\psi$ and $(\bar{\psi},q)$
are inverse equivalences over $V$. But this is clear: $\bar{\psi}\psi(w,z)$
is just $shrink_{0,\half}(w)$, and if $q(w,s)=z$, then
\begin{eqnarray}
\psi(\bar{\psi}(w,s),q(w,s)) & = & \psi(below_s(w),z) \nonumber\\
& = & \left((shrink_{0,\half}(below_s(w))\cup 
shrink_{\half,1}(\lambda(z)),\threeq\right), \nonumber
\end{eqnarray}
and linearly deforming all the segments to their original 
locations, and simultaneously deforming $\threeq$ to $s$ 
linearly, describes a homotopy over $V$ of $\psi\circ(\bar{\psi},q)$ 
to the identity.

The proof of part (2.) rests on the fact that the inclusion 
$F_{j-1}\CYSX\hookrightarrow F_{j}\CYSX$ is a cofibration, which
comes from the fact that $X$ is nondegenerately based. Let $W$
be a neighborhood of the basepoint $*$ in $X$ and let $K_t:X\rar X$
be a based homotopy where $K_0=id$ and $K_1(W)=\{*\}$. Let $L_t$ be
a linear deformation of $[0,1]$ from the identity to the map
$$L_1(t)=\cases{0,&\hbox{if $t\le\qtr$;}\cr
\noalign{\vskip2pt}
2t-\half,&\hbox{if $\qtr\le t\le\threeq$; and}\cr
\noalign{\vskip2pt}
1,&\hbox{if $t\ge\threeq$.}}$$
Use the same symbol $L_t$ to denote the induced homotopy on $S^1$.
Then $J_t=K_t\wedge L_t$ is a deformation of $\Sigma X=X\wedge S^1$ which
collapses a neighborhood $W'=W\wedge\big([0,\qtr)\cup(\threeq,1]\big)$
of the basepoint. Thus let
$$U=\left\{\big[y_i,[x_i,s_i]\big]_{1\le i\le j}\big|
[x_i,s_i]\in W'\hbox{ for some $i$}\right\},$$
and use the functoriality of $C_1(-,-)$ to define $h_t(z)=C(1_Y,J_t)(z)$.
For any $z=\big[y_i,[x_i,s_i]\big]$ in $U$, $J_1([x_i,s_i])$ will
be $*$ for at least one index $i$, and so 
$J_1(z)\in F_{j-1}\CYSX$. It is clear that $J_t$ 
preserves $F_{j-1}\CYSX$ and so part (2a) is complete.

If $(w,s)\in q^{-1}(U)$ and $w=\big[a_i,b_i,y_i,x_i\big]_{1\le i\le j}$, 
define
$$H_t(w,s)= \big(\big[(1-t)a_i+ta_i',(1-t)b_i+tb_i',y_i,K_t(x_i)\big]
_{1\le i\le j},\ s\big),$$
where $a_i'= a_i+\qtr(b_i-a_i)$ and $b_i'= b_i-\qtr(b_i-a_i)$. It is 
straightforward to verify that $qH_t=h_t q$ and so (2b.) is complete.

Finally, the restriction of $H_1$ to fibers fits into a homotopy-commutative 
diagram
$$\xymatrix{
q^{-1}(z) \ar[rr]^{H_1} \ar[d]_{\bar{\psi}} & &q^{-1}(h_1(z)) \ar[d]^{\bar{\psi}} \\
\CbarYX \ar[rr]_{\xi\circ C_1(1_Y,K_1)} & & \CbarYX
}$$
where we have already shown that the maps $\bar{\psi}$ are homotopy 
equivalences, and where $\xi$ is multiplication by the element
$$\big[a'_i,b'_i,y_i,K_1(x_i)\big]_{1\le i \le j\hbox{ and }b'_i\le s<b_i}$$
in the $H$-space structure on $\CbarYX$. Since $\CbarYX$ is connected 
(because $X$ is) this is a homotopy equivalence. This completes the proof 
of (2c.), and hence $q$ is a quasifibration.

More can be said. By extending and iterating the definition and theorem, we 
can prove
\begin{cor}
Let $Y$ be a weak metric space and $X$ a nondegenerately based space.
For each $n\ge 1$ there is a space $C_n(Y,X)$ and a pair of maps 
$$ C(\RnxY,X) {\buildrel {\phi_n} \over \longleftarrow} C_n(Y,X)
{\buildrel {\alpha_n} \over \longrightarrow} \OnCYSnX $$
such that:
\begin{enumerate}
\item $C_n(-,-)$ is functorial with respect to based maps in the second 
variable and injective maps in the first variable, and $\phi_n$ and 
$\alpha_n$ are natural;
\item $\phi_n$ is a homotopy equivalence; and
\item $\alpha_n$ is a weak homotopy equivalence if $X$ is path-connected.
\end{enumerate}
\end{cor}

There is an evident action of the little $n$-cubes operad ${\mathcal C}_n$ of 
\cite{may} on all the spaces appearing in the Corollary, and $\phi_n$ and 
$\alpha_n$ can be seen to be ${\mathcal C}_n$-maps.

It is also true (and proved in \cite{c}) that when $X$ is not path connected, 
$\alpha_n$ is a group-completion for $n\ge 2$.

\end{document}